\documentclass[preprint,12pt]{elsarticle}
\usepackage{amssymb}

\journal{Computer Physics Communications}

\newtheorem{defn}{\emph{Definition} }

\begin{document}

\begin{frontmatter}

\title{A Phase-Fitted Runge-Kutta-Nystr\"{o}m method for the Numerical Solution of
 Initial Value Problems with Oscillating Solutions}

\author[UoP]{D. F. Papadopoulos}
\ead{dimpap@uop.gr}

\author[UoP]{Z. A. Anastassi}
\ead{zackanas@uop.gr}

\author[UoP]{T. E. Simos\fnref{simos}}
\fntext[simos]{Highly Cited Researcher, Active Member of the
European Academy of Sciences and Arts, Address: Dr. T.E. Simos, 26
Menelaou Street, Amfithea - Paleon Faliron, GR-175 64 Athens,
GREECE, Tel: 0030 210 94 20 091}
\ead{tsimos.conf@gmail.com,
tsimos@mail.ariadne-t.gr}

\address[UoP]{Laboratory of Computer Sciences,\\
Department of Computer Science and Technology,\\
Faculty of Sciences and Technology, University of Peloponnese\\
GR-22 100 Tripolis, GREECE}

\begin{abstract}
A new Runge-Kutta-Nystr\"{o}m method, with phase-lag of order
 infinity, for the integration of second-order periodic initial-value
 problems is developed in this paper. The new method is based on the
Dormand and Prince Runge-Kutta-Nystr\"{o}m
 method of algebraic order four\cite{pa}. Numerical illustrations indicate that the new method is
 much more efficient than the classical one.
\end{abstract}

\begin{keyword}
Runge-Kutta-Nystr\"{o}m methods; Phase-fitted; Initial-value
problems; Phase-lag infinity

\PACS 02.60

% MSC codes here, in the form: \MSC code \sep code
% or \MSC[2008] code \sep code (2000 is the default)

\end{keyword}

\end{frontmatter}

\section{Introduction}

    \quad In this paper we study a special Runge-Kutta-Nystr\"{o}m method of
    Dormand $et$ $al.$\cite{pa} for integrating systems of ODEs of the form
    \begin{equation}
    \frac{d^2u(t)}{dt^2} = f(t,u(t))
    \label{eq:no1}
    \end{equation}
    for which it is known in advantage that their solution is periodic
    or oscillating.

    Several authors in their papers (\emph{for example see} [3,7-10]) have developed
    Runge-Kutta-Nystr\"{o}m methods with the purpose of making the
    phase-lag of the method smaller.

    The \emph{phase-lag} of a method, first defined by Brusa and
    Nigro \cite{b} at 1980. Van der Houwen and Sommeijer \cite{c} proposed
    second-order $m$-stage methods (with $m=4,5,6$) and phase-lag
    order $q=6,8,10$ respectively. They also derived some
    third-order methods with phase-lag order $6,8,10$. In \cite{c,e}
    Chawla and Rao have constructed Numerov-type methods with
    minimal phase-lag for the numerical integration of second-order
    initial-value problems. Simos $et$ $al.$ \cite{h} obtain
    fourth-order Runge-Kutta-Nystr\"{o}m with minimal phase-lag of
    order eigth. He also derived in \cite{i} a Runge-Kutta-Fehlberg method of
    order infinity.

    In the present paper and based on the requirements of infinite
    order of phase-lag, we will construct a phase-fitted four-stage
    Runge-Kutta-Nystr\"{o}m which is based on the coefficients of the
    well-known Runge-Kutta-Nystr\"{o}m Dormand $et$ $al.$ \cite{pa} method
    of algebraic order 4.

\section{Phase lag analysis for Runge-Kutta-Nystr\"{o}m methods}

\quad The general m-stage method for the equation
\begin{equation}
\frac{d^2u(t)}{dt^2} = f(t,u(t)) \label{eq:no2}
\end{equation}
is of the form
\begin{eqnarray}
\nonumber &u_n&^{(0)}=u_{n-1},\qquad
u_n^{(i)}=u_{n-1}+h\hat{u}_{n-1}+h^2\sum_{j=1}^{i}{b_jf_j},\\
&u_n&=u_n^{(m)}, \qquad \quad
\hat{u}_n=\hat{u}_{n-1}+h\sum_{j=1}^{i}{\hat{b}_jf_j},
\label{eq:no3}
\end{eqnarray}
where
\begin{equation}
f_i=f(t_{n-1}+c_ih,u_{n-1}+hc_i\hat{u}_{n-1}+h^2\sum_{j=1}^{i-1}{\alpha_{i,j}u_n^{(j)}})
\label{eq:no4}
\end{equation}
and $c_1=0$ and $c_m=1$

The above expressions are presented using the well-known Butcher
table, given below:
\begin{table}
\begin{tabular}{cccccc}
\hline
0\\
$c_2$ & $\alpha_{21}$\\
$c_3$ & $\alpha_{31}$ & $\alpha_{32}$\\
\vdots & \vdots & \vdots\\
$c_m$ & $\alpha_{m,1}$ & $\alpha_{m,2}$ & \ldots & $\alpha_{m,m-1}$\\
\hline
 &$b_1$ & $b_2$ & \ldots & $b_{m-1}$ & $b_m$\\
 &$\hat{b}_1$ & $\hat{b}_2$ & \ldots & $\hat{b}_{m-1}$ & $\hat{b}_m$\\
\end{tabular}
\caption{m-stage Runge-Kutta-Nyst\"{o}m method} \label{tab:a}
\end{table}

In order to develop the new method, we use the test equation,
\begin{equation}
\frac{d^2u(t)}{dt^2}=(iv)^2u(t) \Longrightarrow u''(t)=-v^2u(t),
\qquad v \in R \label{eq:no5}
\end{equation}

By applying the general method (\ref{eq:no3}) to the test equation
(\ref{eq:no5}) we obtain the numerical solution
\begin{equation}
\biggl[\begin{array}{c}
  u_n \\
  h\hat{u}_n
\end{array}\biggl]=D^n\biggl[\begin{array}{c}
  u_0 \\
  h\hat{u}_0
\end{array}\biggl], \quad D=\biggl[\begin{array}{cc}
                                     A(z^2) & B(z^2) \\
                                     A'(z^2) & B'(z^2)
                                   \end{array}
\biggl], \quad z=vh, \label{eq:no6}
\end{equation}
where $A,B,A',B'$ are polynomials in $z^2$, completely determined by
the parameters of method (\ref{eq:no3})

The exact solution of (\ref{eq:no5}) is given by
\begin{equation}
u(t_n)=\sigma_1[exp(iv)]^n+\sigma_2[exp(-iv)]^n, \label{eq:no7}
\end{equation}
where \begin{eqnarray} \nonumber
\sigma_{1,2}=\frac{1}{2}[u_0\pm\frac{(i\hat{u}_0)}{v}] \quad or
\quad \sigma_{1,2}=|\sigma|exp(\pm i\chi).
\end{eqnarray}
Substituting in (\ref{eq:no7}), we have
\begin{equation}
u(t_n)=2|\sigma|cos(\chi+nz). \label{eq:no8}
\end{equation}
Furthermore we assume that the eigenvalues of $D$ are
$\varrho_1,\varrho_2$, and the consequent eigenvectors are
$[1,v_1]^T,[1,v_2]^T$,\\
where $v_i={A'}/{(\rho_i-B')},i=1,2$. The numerical solution of
(\ref{eq:no5}) is
\begin{equation}
u_n=c_1\rho_1^{n}+c_2\rho_2^{n}, \label{eq:no9}
\end{equation}
where
\begin{eqnarray}
\nonumber c_1=-\frac{v_2u_0-h\hat{u}_0}{v_1-v_2},\quad
c_2=-\frac{v_1u_0-h\hat{u}_0}{v_1-v_2}.
\end{eqnarray}

If $\rho_1,\rho_2$ are complex conjugate, then $c_{1,2}=|c|exp(\pm
iw) \; and \; \rho_{1,2}=|\rho|exp(\pm ip)$. By substituting in
(\ref{eq:no9}), we have
\begin{equation}
u_n=2|c||\rho|^ncos(w+np). \label{eq:no10}
\end{equation}
From equations (\ref{eq:no8}) and (\ref{eq:no10}) we take the
following definition.
\begin{defn}
\label{def:no1} (Phase-lag). Apply the RKN method (\ref{eq:no3}) to
the general method (\ref{eq:no5}). Then we define the phase-lag
$\Phi(z)=z-p$. If $\Phi(z)=O(z^{q+1})$, then the RKN method is said
to have phase-lag order $q$.
\end{defn}

In addition, the quantity $a(z)=1-|\rho|$ is called
\emph{amplification error}.

Let us denote
\begin{eqnarray}
\nonumber R(z^2)&=&tr(D)=A(z^2)+B'(z^2)\\
Q(z^2)&=&det(D)=A(z^2)B'(z^2)-A'(z^2)B(z^2) \label{eq:no11}
\end{eqnarray}

where $z=vh$. From Definition \ref{def:no1} it follows that
\begin{equation}
\Phi(z)=z-arcoss\biggl(\frac{R(z^2)}{2\sqrt{Q(z^2)}}\biggl),\qquad
|\rho|=\sqrt{Q(z^2)}. \label{eq:no12}
\end{equation}

We can also put forward an alternative definition for the case of
infinite order of phase lag.
\begin{defn}
\label{def:no2} (Phase-lag of order infinity). To obtain phase-lag
of order infinity the relation
$\Phi(z)=z-arccos\biggl(\frac{R(z^2)}{2\sqrt{Q(z^2)}}\biggl)=0$ must
hold.
\end{defn}

\section{Derivation of the new Runge-Kutta-Nystr\"{o}m method}

\quad In this section we construct a 4-stage explicit
Runge-Kutta-Nystr\"{o}m method (presented in Table 1), based on
$R(z^2)$ and $Q(z^2)$. Now let us rewrite R and Q in the following
form
\begin{eqnarray}
\nonumber R(z^2)=2-r_1z^2+r_2z^4-r_3z^6+\ldots+r_iz^{2i}=0\\
Q(z^2)=1-q_1z^2+q_2z^4-q_3z^6+\ldots+q_iz^{2i}=0 \label{eq:no13}
\end{eqnarray}

By computing the polynomials $A,B,A',B'$ and therefore $R$ and  $Q$
in terms of RKN parameters we obtain the following expressions
\\
\\
$A(z^2)$ = $1+b_{{4}}a_{{4,3}}a_{{3,2}}a_{{2,1}}{z}^{8}+ (
-b_{{4}}a_{
{4,2}}a_{{2,1}}-b_{{3}}a_{{3,2}}a_{{2,1}}-b_{{4}}a_{{4,3}}a_{{3,1}}-b_
{{4}}a_{{4,3}}a_{{3,2}} ) {z}^{6}+( b_{{2}}a_{{2,1}
}+b_{{4}}a_{{4,1}}+b_{{4}}a_{{4,2}}+b_{{3}}a_{{3,1}}+b_{{4}}a_{{4,3}}+
b_{{3}}a_{{3,2}} ) {z}^{4}+ ( -b_{{4}}-b_{{1}}-b_{{3 }}-b_{{2}} )
{z}^{2}$
\\
\\
$B(z^2)$ = $1-b_{{4}}a_{{4,3}}a_{{3,2}}c_{{2}}{z}^{6}+ (
b_{{4}}a_{{4,
3}}c_{{3}}+b_{{4}}a_{{4,2}}c_{{2}}+b_{{3}}a_{{3,2}}c_{{2}} )
{z}^{4}+ ( -b_{{3}}c_{{3}}-b_{{4}}c_{{4}}-b_{{2}}c_{{2}}
 ) {z}^{2}$
\\
\\
$A'(z^2)$ = $\hat{b}_{{4}}a_{{4,3}}a_{{3,2}}a_{{2,1}}{z}^{8}+ (
-\hat{b}_{{3}}a_{{3
,2}}a_{{2,1}}-\hat{b}_{{4}}a_{{4,3}}a_{{3,1}}-\hat{b}_{{4}}a_{{4,3}}a_{{3,2}}-\hat{b}_{{
4}}a_{{4,2}}a_{{2,1}} ) {z}^{6} +  ( \hat{b}_{{2}}a_{{2,1}}+
\hat{b}_{{3}}a_{{3,1}}+\hat{b}_{{3}}a_{{3,2}}+\hat{b}_{{4}}a_{{4,1}}+\hat{b}_{{4}}a_{{4,2}}+\hat{b}_
{{4}}a_{{4,3}} ) {z}^{4}+ (
-\hat{b}_{{4}}-\hat{b}_{{2}}-\hat{b}_{{1}} -\hat{b}_{{3}} ) {z}^{2}$
\\
\\
$B'(z^2)=1$-$\hat{b}_{{4}}a_{{4,3}}a_{{3,2}}c_{{2}}{z}^{6}+ (
\hat{b}_{{4}}a_{{4,
3}}c_{{3}}+\hat{b}_{{4}}a_{{4,2}}c_{{2}}+\hat{b}_{{3}}a_{{3,2}}c_{{2}}
) {z}^{4}+ (
-\hat{b}_{{3}}c_{{3}}-\hat{b}_{{4}}c_{{4}}-\hat{b}_{{2}}c_{{2}}
 ) {z}^{2}$
\\
\\
$R(z^2)=2$+$b_{{4}}a_{{4,3}}a_{{3,2}}a_{{2,1}}{z}^{8}+ ( -b_{{3}}a_{
{3,2}}a_{{2,1}}-b_{{4}}a_{{4,3}}a_{{3,2}}-b_{{4}}a_{{4,2}}a_{{2,1}}-\hat{b}_
{{4}}a_{{4,3}}a_{{3,2}}c_{{2}}-b_{{4}}a_{{4,3}}a_{{3,1}} ) {z}^{6}+
( b_{{2}}a_{{2,1}}+b_{{3}}a_{{3,2}}+b_{{4}}a_{{4,3}}+\hat{b}
_{{3}}a_{{3,2}}c_{{2}}+\hat{b}_{{4}}a_{{4,3}}c_{{3}}+\hat{b}_{{4}}a_{{4,2}}c_{{2}}+
b_{{3}}a_{{3,1}}+b_{{4}}a_{{4,1}}+b_{{4}}a_{{4,2}} ) {z} ^{4}+ (
-b_{{3}}-b_{{2}}-\hat{b}_{{3}}c_{{3}}-\hat{b}_{{4}}c_{{4}}-\hat{b}_{{2}}c_{
{2}}-b_{{4}}-b_{{1}} ) {z}^{2}$
\\
\\
$Q(z^2)=1$+$(-\hat{b}_{{4}}a_{{4,3}}a_{{3,1}}b_{{2}}c_{{2}}-\hat{b}_{{4}}a_{{4,2}}a_{
{2,1}}b_{{3}}c_{{3}}-\hat{b}_{{2}}a_{{2,1}}b_{{4}}a_{{4,3}}c_{{3}}
$-$\hat{b}_{{3}}a_{{3,2}}a_{{2,1}}b_{{4}}c_{{4}}+b_{{3}}a_{{3,1}}\hat{b}_{{4}}a_{{4,2}}c_{{2}}
-\hat{b}_{{3}}a_{{3,1}}b_{{4}}a_{{4,2}}c_{{2}}-\hat{b}_{{4}}a_{{4,3}}a_{{3,2}}a_{{
2,1}}+b_{{4}}a_{{4,2}}a_{{2,1}}\hat{b}_{{3}}c_{{3}}$-$\hat{b}_{{1}}b_{{4}}a_{{4,3}}a
_{{3,2}}c_{{2}}+b_{{4}}a_{{4,1}}\hat{b}_{{3}}a_{{3,2}}c_{{2}}
-\hat{b}_{{4}}a_{{4,1}}b_{{3}}a_{{3,2}}c_{{2}}+b_{{4}}a_{{4,3}}a_{{3,2}}a_{{2,1}}+b_{{1}}\hat{b}_
{{4}}a_{{4,3}}a_{{3,2}}c_{{2}}$+$b_{{3}}a_{{3,2}}a_{{2,1}}\hat{b}_{{4}}c_{{4}}
+b_{{4}}a_{{4,3}}a_{{3,1}}\hat{b}_{{2}}c_{{2}}+b_{{2}}a_{{2,1}}\hat{b}_{{4}}a_{{4,
3}}c_{{3}}) {z}^{8}+(-b_{{4}}a_{{4,3}}a_{{3,1}}-b
_{{3}}a_{{3,2}}a_{{2,1}}-b_{{4}}a_{{4,2}}a_{{2,1}}$-$b_{{4}}a_{{4,3}}a_{
{3,2}}-b_{{1}}\hat{b}_{{4}}a_{{4,2}}c_{{2}}-b_{{1}}\hat{b}_{{3}}a_{{3,2}}c_{{2}}-b
_{{3}}\hat{b}_{{4}}a_{{4,2}}c_{{2}}+\hat{b}_{{2}}b_{{4}}a_{{4,3}}c_{{3}}-b_{{2}}a_
{{2,1}}\hat{b}_{{3}}c_{{3}}-b_{{2}}a_{{2,1}}\hat{b}_{{4}}c_{{4}}$-$b_{{4}}a_{{4,1}}\hat{b}
_{{3}}c_{{3}}-b_{{4}}a_{{4,1}}\hat{b}_{{2}}c_{{2}}-b_{{4}}a_{{4,2}}\hat{b}_{{3}}c_
{{3}}+\hat{b}_{{3}}a_{{3,2}}b_{{4}}c_{{4}}+\hat{b}_{{4}}a_{{4,1}}b_{{3}}c_{{3}}+\hat{b}_
{{4}}a_{{4,1}}b_{{2}}c_{{2}}+\hat{b}_{{4}}a_{{4,2}}b_{{3}}c_{{3}}$-$b_{{2}}\hat{b}_{
{4}}a_{{4,3}}c_{{3}}+\hat{b}_{{2}}a_{{2,1}}b_{{4}}c_{{4}}+\hat{b}_{{3}}a_{{3,1}}b_
{{4}}c_{{4}}+\hat{b}_{{3}}a_{{3,1}}b_{{2}}c_{{2}}-b_{{3}}a_{{3,1}}\hat{b}_{{4}}c_{
{4}}-b_{{3}}a_{{3,1}}\hat{b}_{{2}}c_{{2}}-b_{{4}}a_{{4,3}}\hat{b}_{{2}}c_{{2}}$-$b_{
{3}}a_{{3,2}}\hat{b}_{{4}}c_{{4}}-b_{{4}}\hat{b}_{{3}}a_{{3,2}}c_{{2}}-b_{{1}}\hat{b}_{{
4}}a_{{4,3}}c_{{3}}+\hat{b}_{{4}}b_{{3}}a_{{3,2}}c_{{2}}+\hat{b}_{{1}}b_{{4}}a_{{4
,3}}c_{{3}}+\hat{b}_{{4}}a_{{4,3}}b_{{2}}c_{{2}}+\hat{b}_{{1}}b_{{4}}a_{{4,2}}c_{{
2}}$+$\hat{b}_{{1}}b_{{3}}a_{{3,2}}c_{{2}}+\hat{b}_{{3}}b_{{4}}a_{{4,2}}c_{{2}}+\hat{b}_{{
2}}a_{{2,1}}b_{{3}}c_{{3}}+\hat{b}_{{3}}a_{{3,2}}a_{{2,1}}+\hat{b}_{{4}}a_{{4,3}}a
_{{3,1}}+\hat{b}_{{4}}a_{{4,2}}a_{{2,1}}+\hat{b}_{{4}}a_{{4,3}}a_{{3,2}}-\hat{b}_{{4}}a_
{{4,3}}a_{{3,2}}c_{{2}})\\
 {z}^{6}+( -\hat{b}_{{4}}b_{{3}}
c_{{3}}+b_{{4}}\hat{b}_{{2}}c_{{2}}-\hat{b}_{{2}}b_{{4}}c_{{4}}-\hat{b}_{{3}}b_{{4}}c_{{
4}}+b_{{2}}\hat{b}_{{4}}c_{{4}}+b_{{3}}\hat{b}_{{2}}c_{{2}}+b_{{3}}\hat{b}_{{4}}c_{{4}}-
\hat{b}_{{1}}b_{{3}}c_{{3}}$+$b_{{4}}\hat{b}_{{3}}c_{{3}}+b_{{1}}\hat{b}_{{3}}c_{{3}}-\hat{b}_{{
4}}b_{{2}}c_{{2}}-\hat{b}_{{1}}b_{{2}}c_{{2}}-\hat{b}_{{1}}b_{{4}}c_{{4}}+b_{{1}}\hat{b}
_{{2}}c_{{2}}+b_{{2}}\hat{b}_{{3}}c_{{3}}+b_{{1}}\hat{b}_{{4}}c_{{4}}-\hat{b}_{{3}}b_{{2
}}c_{{2}}$-$\hat{b}_{{2}}b_{{3}}c_{{3}}-\hat{b}_{{2}}a_{{2,1}}-\hat{b}_{{3}}a_{{3,1}}-\hat{b}_{{
3}}a_{{3,2}}-\hat{b}_{{4}}a_{{4,1}}-\hat{b}_{{4}}a_{{4,2}}-\hat{b}_{{4}}a_{{4,3}}+b_{{2}
}a_{{2,1}}+b_{{4}}a_{{4,1}}+b_{{4}}a_{{4,2}}+b_{{3}}a_{{3,1}}$+$b_{{4}}a
_{{4,3}}+b_{{3}}a_{{3,2}}+\hat{b}_{{3}}a_{{3,2}}c_{{2}}+\hat{b}_{{4}}a_{{4,3}}c_{{
3}}+\hat{b}_{{4}}a_{{4,2}}c_{{2}}) {z}^{4}+ ( -b_{{2}}-\hat{b}_
{{4}}c_{{4}}+\hat{b}_{{2}}-b_{{4}}+\hat{b}_{{1}}$+$\hat{b}_{{3}}-\hat{b}_{{2}}c_{{2}}-b_{{1}}-b_
{{3}}-\hat{b}_{{3}}c_{{3}}+\hat{b}_{{4}} ) {z}^{2}$
\\
\\
where $z=\nu h$\\

As it has already been defined, in order to have phase-lag of order
infinity, the following relation must hold:\\
\begin{equation}
\Phi(z)=z-arccos\biggl(\frac{R(z^2)}{2\sqrt{Q(z)^2}}\biggl)=0
\label{eq:no14}
\end{equation}

By applying $R(z^2)$ and $Q(z^2)$ to the formula of the direct
calculation of the phase lag (\ref{eq:no12}) and substituting the
following coefficients that have been used by Dormand et al. in \cite{pa} :\\

\begin{eqnarray}
\nonumber &\alpha_{21}&=\frac{1}{32}, \qquad
\alpha_{31}=\frac{7}{1000}, \qquad \alpha_{32}=\frac{119}{500},
\qquad \alpha_{41}=\frac{1}{14}, \qquad \alpha_{42}=\frac{8}{27},
\\
\nonumber &c_2&=\frac{1}{4}, \qquad c_3=\frac{7}{10}, \qquad c_4=1,
\\
\nonumber &b_1&=\frac{1}{14}, \qquad b_2=\frac{8}{27}, \qquad
b_3=\frac{25}{189}, \qquad
b_4=0, \\
\nonumber &\hat{b}_1&=\frac{1}{14}, \qquad \hat{b}_2=\frac{32}{81},
\qquad \hat{b}_3=\frac{250}{567}, \qquad \hat{b}_4=\frac{5}{54},
\end{eqnarray}
After satisfying relation (\ref{eq:no14}), we have:
\begin{eqnarray}
\nonumber \Phi(z)&=&z-arcoss\biggl(\frac{R(z^2)}{2\sqrt{Q(z)^2}}\biggl)=0\Rightarrow \\
\nonumber a_{4,3}&=&-\frac
{5}{5292}\frac{1}{289z^4-6800z^2+40000z^4}(54621\,{z}^{8}-4793320\,{z}^{6}+99172960\,{z}^{4}\\
\nonumber &+&5179680\,{z}^{4}
 \left( \sin \left( z \right)  \right) ^{2}-768268800\,{z}^{2}+4043520 \,{z}^{2} \left( \sin \left( z \right)
\right) ^{2}\\
\nonumber &+&1866240000- 559872000\, \left( \sin \left( z \right)
\right)
^{2}+24\,(-654383577600\,{z}^{6}\\
\nonumber &+&212348252160000\,{z}^{4}-1366377865200\,{z}^{8}-
1710031785\,{z}^{12}\\
\nonumber &+&89285428680\,{z}^{10}-202307339750400\,{z}^{4}
 \left( \sin \left( z \right)  \right) ^{2}\\
 \nonumber &+&2023399802880000\,{z}^{2}
 \left( \sin \left( z \right)  \right) ^{2}-2015539200000000\,{z}^{2}\\
 \nonumber &+& 581660870400\,{z}^{6} \left( \sin \left(
z \right)  \right) ^{2}+ 1319799592800\,{z}^{8} \left(
\sin \left( z \right)  \right) ^{2}\\
\nonumber &+&1710031785\,{z}^{12} \left( \sin \left( z \right)
\right) ^{2}- 89285428680\,{z}^{10} \left( \sin \left( z \right)
\right) ^{2}\\
\nonumber &+& 46578272400\,{z}^{8} \left( \sin \left( z \right)
\right) ^{4}+72722707200\,{z}^{6} \left( \sin \left( z \right)
\right) ^{4}\\
\nonumber &-& 10040912409600\,{z}^{4} \left( \sin \left( z \right)
\right) ^{4}+ 544195584000000\, \left( \sin \left( z \right) \right)
^{4}\\
\nonumber &-&7860602880000\,{z}^{2} \left( \sin \left( z \right)
\right) ^{4}+
6046617600000000\\
&-&6590813184000000\, \left( \sin \left( z \right)
 \right) ^{2})^{1/2}
) \label{eq:no15}
\end{eqnarray}
\\
The Taylor expansion series for $a_{4,3}$, which is given from the
above formula  is :
\begin{eqnarray}
\nonumber a_{4,3}&=&{\frac {25}{189}}-{\frac
{43}{2400}}\,{z}^{2}-{\frac {1531}{30240000}} \,{z}^{4}-{\frac
{3273029}{36288000000}}\,{z}^{6}\\
&+&{\frac {59772887431} {9699782400000000}}\,{z}^{8} +\cdots .
\label{eq:no16}
\end{eqnarray}

\section{Numerical examples}
\quad In this section we will apply our method to three problems. We
are going to compare our results with those derived by using the
high order method of embedded Runge-Kutta-Nystr\"{o}m $4(3)4$ method
of \emph{Dormand and Prince (see \cite{pa})}.\\

One way to measure the efficiency of the method is to compute the
accuracy in the decimal digits, that is $-log_{10}$(\emph{maximum error through the integration intervals})
\\
\\
$acc(T)=-log_{10}(max|u(t_n)-u_n|),\quad where \quad t_n=1+nh,\quad
n=1,2,\ldots,\frac{T-1}{h}$ and $u(t)$ is the vector of the solution.
\\
\\
{Table 2} shows the accuracy for the two methods. In our
computations we have two step values, for Problems 1 and 2,
$h=0.025$ and $h=0.050$, and for Problems 3 and 4, $h=0.25$ and
$h=0.50$.
\\

\emph{\textbf{Problem 1.}(Inhomogeneous equation)}
\begin{eqnarray}
\nonumber \frac{d^2u(t)}{dt^2}=-\nu^2u(t)+(\nu^2-1)sin(t),\qquad
u(0)=1, \quad u'(0)=\nu+1,
\end{eqnarray}
where $t\geq 0$ and $\nu=10$.\\
The analytical solution is $u(t)=cos(\nu t)+sin(\nu t)+sin(t)$
\\

\emph{\textbf{Problem 2.}(Two-Body problem)}
\begin{eqnarray}
\nonumber u''=-\frac{u}{(u^2+z^2)^{3/2}},\quad
z''=-\frac{z}{(u^2+z^2)^{3/2}}
\end{eqnarray}
where $\qquad u(0)=1,\quad u'(0)=0,\quad z(0)=0,\quad z'(0)=1 \;\; and \;\; \nu=1$

The analytical solution is $u(t)=cos(t) \quad and \quad z(t)=sin(t)$
\\

\emph{\textbf{Problem 3.}(Duffing equation)}
\begin{eqnarray}
\nonumber \frac{d^2u(t)}{dt^2}=-u(t)-(u(t))^3+Bcos(\nu t)
\end{eqnarray}
where $B=0.002$ and $\nu=1.01$.\\
The analytical solution is $u(t)=A_1cos(\nu t)+A_3cos(3\nu
t)+A_5cos(5\nu
t)+A_7cos(7\nu t)+A_9cos(9\nu t)$\\
where $A_1=0.200179477536$, $A_3=0.000246946143$, $A_5=0.000000304014$, $A_7=0.000000000374$ and $A_9=0.000000000000$\\

\emph{\textbf{Problem 4.}(Franco and Palacios problem)}
\begin{eqnarray}
\nonumber \frac{d^2u(t)}{dt^2}=-u(t)+\epsilon exp(it),\qquad u(t)\in
C \quad u(0)=1, \quad u'(0)=(1-\frac{1}{2}\epsilon)i,
\end{eqnarray}
where $\epsilon=0.001$ and $\nu=1$\\
The analytical solution is $u(t)=cos(t)+\frac{1}{2}\epsilon
tsin(t)+i[sin(t)-\frac{1}{2}\epsilon tcos(t)]$
\\

\begin{table}
\qquad \qquad \qquad \quad \qquad Our method \qquad \quad \quad Dormand and Prince method \\
\begin{tabular}{lcccccc}
\hline
 & T=100\quad & T=1000\quad & T=5000 & T=100\quad & T=1000\quad & T=5000\\
\hline
\underline{\emph{Problem 1}}\\
h=0.025 & \quad4.2 & \qquad3.2 & 2.5 & \quad2.3 & \qquad1.3 & 0.6\\
h=0.050 & \quad2.7 & \qquad1.7 & 1.0 & \quad1.1 & \qquad0.2 & -0.3\\
\hline
\underline{\emph{Problem 2}}\\
h=0.025 & \quad7.3 & \qquad5.9 & 4.6 & \quad6.5 & \qquad5.1 & 3.8\\
h=0.050 & \quad6.0 & \qquad4.4 & 3.1 & \quad5.2 & \qquad3.6 & 2.3\\
\hline
\underline{\emph{Problem 3}}\\
h=0.25 & \quad5.7 & \qquad5.4 & 5.4 & \quad4.2 & \qquad4.1 & 4.1\\
h=0.50 & \quad4.2 & \qquad3.9 & 3.9 & \quad2.9 & \qquad2.8 & 2.8\\
\hline
\underline{\emph{Problem 4}}\\
h=0.25 & \quad5.2 & \qquad4.3 & 3.4 & \quad3.5 & \qquad2.5 & 1.6\\
h=0.50 & \quad3.8 & \qquad2.8 & 1.9 & \quad2.3 & \qquad1.8 & 0.4\\
\hline

\end{tabular}
\caption{Accuracy for the maximum absolute error for problems 1-4}
\label{tab:b}
\end{table}

\section{Conclusion}
\quad A new fourth order Runge-Kutta-Nystr\"{o}m method with
phase-lag of order infinity is developed in the present paper. The
new method is based on the very well known classical Dormand and
Prince fourth algebraic order Runge-Kutta-Nyst\"{o}m method. The
numerical results show that the new method is much more efficient
for integrating second-order equations with periodic oscillating behavior than the classical one.
\\
\\
\\

\end{document}